\documentclass[12pt,leqno]{article}
\usepackage{amsmath,amssymb,amscd,latexsym}
\newcommand{\hc}{\hat{c}}
\newcommand{\hC}{\hat{C}}
\newcommand{\wCH}{\widehat{CH}}
\newcommand{\wdeg}{\widehat{\deg}}
\newcommand{\tc}{\tilde{c}}
\newcommand{\oD}{\overline{D}}
\newcommand{\oE}{\overline{E}}
\newcommand{\oL}{\overline{L}}
\newcommand{\oG}{\overline{G}}
\newcommand{\Rst}{{\mathbb{R}}}
\renewcommand{\H}{\mathrm{H}}
\newcommand{\Ext}{\mathrm{Ext}}
\usepackage[all]{xy}

\newcommand{\C}{{\mathbb{C}}}

\newcommand{\Q}{{\mathbb{Q}}}
\newcommand{\R}{{\mathbb{R}}}
\newcommand{\rk}{\mathrm{rk}}
\newcommand{\Spec}{\mathrm{Spec}\,}

\newcommand{\Ih}{{\mathcal I}}
\newcommand{\Eh}{{\mathcal E}}

\newcommand{\Lh}{{\mathcal L}}
\newcommand{\Oh}{{\mathcal O}}

\newcommand{\tei}{\, | \,}

\newtheorem{theorem}{Theorem}
\newtheorem{lemma}[theorem]{Lemma}

\newtheorem{remark}[theorem]{Remark}

\newenvironment{proofof}{\noindent {\bf Proof of}}{\mbox{}\hfill$\Box$}
\newenvironment{proof}{\noindent {\bf Proof}}{\mbox{}\hfill$\Box$}
\begin{document}
\title{A quantitative sharpening of Moriwaki's arithmetic Bogomolov inequality}
\author{N. Naumann}
\date{\ }
\maketitle

A. Moriwaki proved the following arithmetic analogue of the Bogomolov unstability theorem. If a torsion-free hermitian coherent sheaf on an arithmetic surface has negative
discriminant then it admits an arithmetically 
destabilising subsheaf. In the geometric situation
it is known that such a subsheaf can be found subject to an additional numerical constraint
and here we prove the arithmetic analogue. We then apply this result to
slightly simplify a part of C. Soul\'e's proof of a vanishing theorem on 
arithmetic surfaces.

\section{Introduction and statement of result}

Let $K$ be a number field with ring of integers ${\mathcal O}_K$ and $X / \Spec ({\mathcal O}_K)$ an arithmetic surface, i.e. a regular, integral, purely two-dimensional scheme, proper and flat over $\Spec ({\mathcal O}_K)$ and with smooth and geometrically connected generic fibre. Attached to a hermitian coherent sheaf on $X$ are the usual characteristic classes with values in the arithmetic Chow-groups $\wCH^i (X)$ (cf. \cite{GS1}, 2.5), and in particular the discriminant of $\oE$
\[
\Delta (\oE) := (1 - r) \hc_1 (\oE)^2 + 2r \hc_2 (\oE) \in \wCH^2 (X)
\]
where $r := \rk (E)$. The arithmetic degree map
\[ \wdeg: \wCH^2(X)_{\Rst}\longrightarrow \Rst \]
is an isomorphism \cite{GS2} and we will use the same symbol to
to denote an element in $\wCH^2(X)_{\Rst}$ and its
arithmetic degree in $\Rst$, see \cite{GS2}, 1.1 for the 
definition of arithmetic Chow-groups with real coefficients $\wCH^*(X)_{\R}$. Following \cite{Mo2} we
define the positive cone of $X$ to be

\[ \hC_{++}(X) := \{ x \in \wCH^1 (X)_{\Rst} \tei x^2 > 0 \; \mbox{and} \; \deg_K (x) > 0 \} \; . \]

Given a torsion-free hermitian coherent sheaf $\oE$
of rank $r\ge 1$ on $X$ and a subsheaf $E'\subseteq E$ we endow $E'$ with the metric induced from $\oE$
and consider the difference of slopes
\[ \xi_{\oE',\oE} := \frac{\hc_1 (\oE')}{\rk (E')} - \frac{\hc_1 (\oE)}{r} \in \wCH^1 (X)_{\R}. \]

Recall that a subsheaf $E'\subseteq E$ is {\em saturated} if the quotient
$E/E'$ is torsion-free. Our main result is the following.

\begin{theorem} \label{theorem}
Let $\oE$ be a torsion-free hermitian coherent sheaf of rank $r\ge 2$ on the arithmetic surface $X$, satisfying
\[
\Delta (\oE) < 0 \; .
\]
Then there is a non-zero saturated subsheaf $\oE' \subseteq \oE$ such that $\xi_{\oE' , \oE} \in \hC_{++} (X)$ and
\begin{equation}\label{ineq}
\xi^2_{\oE',\oE} \ge \frac{-\Delta}{r^2 (r-1)} \; .
\end{equation}

\end{theorem}

\begin{remark} 
The existence of an $\oE'\subseteq \oE$ with $\xi_{\oE' , \oE} \in \hC_{++} (X)$ is the main result
of \cite{Mo2} and means that $\oE'\subseteq \oE$
is arithmetically destabilising with respect to any polarisation of $X$, c.f. {\em loc. cit.}
for more details on this. The new contribution here is the inequality (\ref{ineq}) which is the exact arithmetic analogue of a known geometric result, c.f. for example \cite{HL}, Theorem 7.3.4.
\end{remark}
\begin{remark}
A special case of Theorem \ref{theorem} appears
in disguised form in the proof of \cite{So}, Theorem 2:
Given a sufficiently positive hermitian line bundle $\oL$ on the arithmetic surface $X$
and some non-torsion element $e \in \H ^1(X,L^{-1})\simeq \Ext ^1(L,\Oh_X)$, C. Soul\'e establishes a lower
bound for 
\[ ||e||^2:=\sup_{\sigma: K\hookrightarrow \C} \; ||\sigma(e)||^2_{L^2} \]

by considering the extension determined by $e$

\[ \Eh:  0 \longrightarrow \overline{\Oh_X} \longrightarrow \oE \longrightarrow \oL \longrightarrow 0 \]

and suitably metrised as to have $\hc_1(\oE)=
\oL$ and $2\hc_2(\oE)=\sum_{\sigma} ||\sigma(e)||^2_{L^2}$, hence $\Delta(\oE)=-\oL^2+2 \sum_{\sigma} ||\sigma(e)||^2_{L^2}$ (where we write
$\oL=\hc_1(\oL)$ following the notation of {\em loc. cit.}).\\
If $E_{\overline{\Q}}$ is semi-stable the arithmetic Bogomolov
inequality concludes the proof. Otherwise,
the main point is to show the existence of
of an arithmetic divisor $\oD$ satisfying

\begin{eqnarray}
\deg_K(\oD) & \le & \deg_K(\oL)/2 \mbox{ and}\label{eins}\\
2(\oL-\oD)\oD & \leq & [K:\Q]\cdot ||e||^2, \label{zwei}
\end{eqnarray}

c.f. $(28)$ and $(32)$ of {\em loc. cit.} where
these inequalities are established by some direct argument.
We wish to point out that the existence of some $\oD$ satisfying
(\ref{eins}) and (\ref{zwei}) is a special case of Theorem \ref{theorem}. In fact,
let $\oE'\subseteq \oE$ be as in Theorem \ref{theorem}
and define $\oD:=\oL-\hc_1(\overline{E'})$. We
then compute 
\[ \xi_{\oE',\oE}=\frac{\oL}{2}-\oD \]
and $\xi_{\oE',\oE}\in \hC_{++}(X)$ implies (\ref{eins}). Furthermore, the inequality (\ref{ineq}) in the present case reads

\[ \xi_{\oE',\oE}^2=\frac{\oL^2}{4}+\oD^2-\oL\;\oD\ge\frac{-\Delta}{4}=\frac{\oL^2}{4}-\frac{1}{2}
\sum_{\sigma} ||\sigma(e)||^2_{L^2}\mbox{  , i.e.}\]

\[ 2(\oL-\oD)\oD\leq \sum_{\sigma} ||\sigma(e)||^2_{L^2},\]

hence the trivial estimate $[K:\Q]\cdot ||e||^2\geq
\sum_{\sigma} ||\sigma(e)||^2_{L^2}$ gives (\ref{zwei}).
\end{remark}

 I would like to thank K. K\"unnemann for useful
conversations about a preliminary draft of the present note.

\section{Proof of Theorem \ref{theorem}}

We collect some lemmas first. We call a short exact sequence 
\[ \Eh : 0 \longrightarrow \oE' \longrightarrow \oE \longrightarrow \oE'' \longrightarrow 0
\]
of hermitian coherent sheaves on
$X$ {\em isometric} if the metrics on $E'$ and $E''$ are induced from the one on $E$. This implies that $\hc_1(\oE)=\hc_1(\oE')+\hc_1(\oE'')$ (i.e. $\tc_1(\Eh)=0$). We also have
\[
\hc_2 (\oE) = \hc_2 (\oE' \oplus \oE'') - a (\tilde{c}_2 (\Eh)) \quad \mbox{in} \; \wCH^2 (X) \; ,
\]

where
\[
a : \tilde{A}^{1,1} (X_{\Rst}) \longrightarrow \wCH^2 (X)
\]
is the usual map \cite{SABK}, chapter III.

\begin{lemma} \label{1}
  If 
\[
\Eh : 0 \longrightarrow \oE' \longrightarrow \oE \longrightarrow \oE'' \longrightarrow 0
\]
is an isometric short exact sequence of hermitian coherent sheaves on $X$ with ranks $r' , r , r''\ge 1$ and discriminants $\Delta' , \Delta , \Delta''$, then
\[
\frac{\Delta'}{r'} + \frac{\Delta''}{r''} - \frac{\Delta}{r} = \frac{rr'}{r''} \xi^2_{\oE',\oE} + 2a (\tc_2 (\Eh)) \quad \mbox{in} \; \wCH^2 (X)_{\R} \; .
\]
\end{lemma}

\begin{proof} We omit the computation using the 
formulas for $\hc_i(\oE)$ recalled above which
shows that the left hand side of the stated equality 
equals

\[ \hc_1(\oE)^2\left( \frac{r-1}{r}+\frac{1-r'}{r'}\right)+\hc_1(\oE'')^2\left( \frac{r-1}{r}+\frac{1-r''}{r''}\right)+ \]
\[ + \hc_1(\oE')\hc_1(\oE'')\left( \frac{2(r-1)}{r}-2 \right) + 2a(\tc(\Eh)).\]

Similarly one writes $\xi_{\oE',\oE}^2$ as a rational linear combination of $\hc_1(\oE)^2,\hc_1(\oE'')^2$ and $\hc_1(\oE')\hc_1(\oE'')$ and comparing 
the results, the stated formula drops out.
\end{proof}

\begin{lemma}
\label{2}
  For $\Eh$ as in Lemma \ref{1} and $\oG'' \subseteq \oE''$ a saturated subsheaf of rank $s\ge 1$ carrying the induced metric, put 
\[
\oG := \ker (E \longrightarrow E'' \longrightarrow E'' / G'') \subseteq \oE
\]
with the induced metric. Then
\[
\xi_{\oG,\oE} = \frac{r' (r'' - s)}{(r' + s) r''} \xi_{\oE',\oE} + \frac{s}{r'+s} \xi_{\oG'' , \oE''} \quad \mbox{in} \; \wCH^1 (X)_{\R} \; .
\]
\end{lemma}

Observe that the coefficients in the last expression are non-negative rational numbers. 

\begin{proof} We have a commutative diagram
with exact rows and columns

\[ \xymatrix{ & & 0  & 0 & \\
 & & \overline{E/G} \ar[u]\ar^{\simeq}[r] & \overline{E''/G''} \ar[u] \\
\Eh:0 \ar[r] & \oE' \ar[r] & \oE \ar[u]\ar[r] & \oE'' \ar[r] \ar[u] & 0 \\
0 \ar[r] & \overline{H} \ar^{\simeq}[u] \ar[r] & \oG \ar[r]\ar[u] & \oG''\ar[u]\ar[r] & 0 \\
 & & 0 \ar[u] & 0. \ar[u] & } \]

Here, we have endowed $E/G,E''/G''$ and $H$ with the metrics induced from
$\oE$,$\oE''$ and $\overline{G}$, hence all rows and columns are isometric by definition. A minor point
to note is that with this choice of metrics the two indicated isomorphisms are isometric, indeed this only means that taking sub- (resp. quotient-)metrics is transitive. One has 
\[ \xi_{\oE',\oE}=\frac{r''\hc_1(\oE')-r'\hc_1(\oE'')}{r'r} \]
and analogously for any isometric exact sequence in place of $\Eh$.
Using this and the diagram one writes
both sides
of the stated equality as a $\Q$-linear combination of
$\hc_1(\oE'),\hc_1(\oG'')$ and $\hc_1(\overline{E''/G''})$ to obtain the same result, namely

\[ \frac{r''-s}{(r'+s)r}\hc_1(\oE')+ 
\frac{r''-s}{(r'+s)r}\hc_1(\oG'')-\frac{1}{r}\hc_1(\overline{E''/G''}).\]

\end{proof}

Finally, we will need the following observation about the intersection theory
on $X$ where, for $x \in \hC_{++} (X)$, we write $|x| := (x^2)^{1/2} \in \Rst^+$.

\begin{lemma}\label{cone}

The subset $\hC_{++}(X)\subseteq\wCH^1(X)_{\Rst}$ is
an open cone, i.e. $x,y\in \hC_{++}(X)$ and 
$\lambda\in\Rst^+$ implies that $x+y,\lambda
x\in\hC_{++}(X)$. For $x,y\in\hC_{++}$ we have 
$|x+y|\ge |x|+|y|$.
\end{lemma}

\begin{proof} This is \cite{Mo2}, (1.1.2.2) except
for the final assertion which is obvious if
$x\in\Rst y$ and we can thus assume that $V:=\Rst x+\Rst y\subseteq\wCH^1(X)_{\Rst}$ is two-dimensional.
We claim that the restriction of the intersection-pairing
makes $V$ a real quadratic space of type $(1,-1)$.
As we have $x\in V$ and $x^2>0$ we only have to exhibit some $v\in V$
with $v^2<0$. To achieve this let $h\in\wCH^1(X)_{\Rst}$ be the first arithmetic Chern class 
of some sufficiently positive hermitian line bundle on $X$ such that the arithmetic Hodge index 
theorem holds for the Lefschetz operator defined by $h$, c.f. \cite{GS2}, Theorem 2.1, ii). Then $a:=xh$ (resp. $b:=yh$) are non-zero real 
numbers for otherwise we would have $x^2<0$ (resp. $y^2<0$). Thus $v:=\frac{x}{a}-\frac{y}{b}\in V$ satisfies $v\neq 0$ and $vh=0$ , hence 
$v^2<0$.\\
Fix a basis $e,f\in V$ with $e^2=1, f^2=-1$
and write
\[ x=\alpha e + \beta f\mbox{ and}\]
\[ y=\gamma e+\delta f.\]

To show that $|x+y|\ge |x|+|y|$ we can assume,
changing both the signs of $x$ and $y$ if necessary, that $\alpha>0$. We then claim that $\gamma>0$.
For otherwise there would be $\lambda_1,\lambda_2\in\Rst^+$
such that $v:=\lambda_1x+\lambda_2 y$ would
have $e$- coordinate equal to zero, hence $v^2\leq 0$ contradicting the fact that either $-v$ or $v$ lies in $\hC_{++}(X)$ (depending on whether or not
we changed the signs of $x$ and $y$ above).\\
From $x^2=\alpha^2-\beta^2$, $y^2=\gamma^2-\delta^2>0$ we obtain $\alpha=|\alpha|\ge |\beta|$
and $\gamma=|\gamma|\ge|\delta|$ and then $\alpha\gamma\ge |\beta\delta| \ge \beta\delta$, i.e.

\begin{equation}\label{sternchen}
xy=\alpha\gamma-\beta\delta\ge 0.
\end{equation}

To conclude, we use the following chain of equivalent statements

\[ |x+y|\ge |x|+|y| \Leftrightarrow \]
\[ (x+y)^2-(|x|+|y|)^2\ge 0 \Leftrightarrow \]
\[ 2xy-2|x||y|\ge 0 \Leftrightarrow \]
\[ xy\ge |x||y| \stackrel{(\ref{sternchen})}{\Leftrightarrow} \]
\[ (xy)^2\ge |x|^2|y|^2 \Leftrightarrow \]
\[ (\alpha\gamma-\beta\delta)^2\ge (\alpha^2-\beta^2)(\gamma^2-\delta^2) \Leftrightarrow \]
\[ \alpha^2\gamma^2+\beta^2\delta^2-2\alpha\beta\gamma\delta\geq \alpha^2\gamma^2-\alpha^2\delta^2-\beta^2\gamma^2+\beta^2\delta^2 \Leftrightarrow \]
\[ 2\alpha\beta\gamma\delta\leq \alpha^2\delta^2+\beta^2\gamma^2 \Leftrightarrow \]
\[ 0\leq (\alpha\delta-\beta\gamma)^2. \]

\end{proof}

\begin{proofof} Theorem \ref{theorem}.
We first remark that for a torsion-free hermitian coherent sheaf $\overline{F}$ of rank one on $X$ 
we always have $\Delta(\overline{F})\ge 0$. 
In fact,
\[
F \simeq \Lh \otimes \Ih_Z
\]
for some line-bundle $\Lh$ and $\Ih_Z$ the ideal sheaf of some closed subscheme $Z \subseteq X$ of codimension $2$. This becomes an isometry for the trivial metric on $\Ih_Z$ and a suitable metric on $\Lh$ (since $\Ih_Z$ is trivial on the generic fibre of $X$). Then
\[
\Delta (\overline{F}) = 2 \hc_2 (\overline{\Lh} \otimes \Ih_Z) = 2 \hc_2 (\Ih_Z) = 2 \, \mbox{length} (Z) \ge 0 \; .
\]
By the main result of \cite{Mo2}, there is $0 \neq \oE' \subseteq \oE$ saturated such that $\xi_{\oE',\oE} \in \hC_{++} (X)$. 
We can assume that, as $E'$ varies through these subsheaves, the real numbers $\xi^2_{\oE',\oE}$ remain bounded for otherwise there is nothing to prove. So we can choose $0 \neq \oE' \subseteq \oE$ saturated with $\xi_{\oE' , \oE} \in \hC_{++} (X)$ and $\xi^2_{\oE',\oE}$ maximal subject to these conditions.
Put $E'' := E / E'$ and consider the isometric exact sequence
\[
\Eh : 0 \longrightarrow \oE' \longrightarrow \oE \longrightarrow \oE'' \longrightarrow 0
\]
with discriminants $\Delta' , \Delta , \Delta''$ and ranks $r' , r, r''$. We claim that 
$\Delta'\ge 0$. This is clear in case $r=2$ 
from the remark made at the beginning of the proof. In case  $r\ge 3$ we assume that $\Delta'<0$ and we let $\overline{G} \subseteq \oE'$ be a saturated subsheaf with $\xi_{\oG,\oE'} \in \hC_{++}$.
Then $\overline{G} \subseteq \oE$ is saturated and using lemma \ref{cone} 
we get
\[
|\xi_{\oG,\oE}| = |\xi_{\oG,\oE'} + \xi_{\oE',\oE}| \ge |\xi_{\oG,\oE'}| + |\xi_{\oE',\oE}| > |\xi_{\oE',\oE}| 
\]
contradicting the maximality of $|\xi_{\oE',\oE}|$. So we have indeed $\Delta' \ge 0$.
Assume now, contrary to our assertion, that
\begin{equation}
  \label{eq:1}
  \frac{\Delta}{r} < -r (r-1) \xi^2_{\oE',\oE} \; .
\end{equation}
Then from Lemma \ref{1}, $\Delta' \ge 0$, (\ref{eq:1}) and $\tc_2 (\Eh) \le 0$ (\cite{Mo1}, 7.2) we get
\begin{gather}
  \frac{\Delta''}{r''} \le \frac{\Delta}{r} + \frac{rr'}{r''} \xi^2_{\oE',\oE} < \left( -r (r-1) + \frac{rr'}{r''} \right) \xi^2_{\oE',\oE} \notag \\
= - r^2 \frac{r''-1}{r''} \xi^2_{\oE',\oE} \le 0 \; , \notag
\end{gather}
hence $\Delta'' < 0$. By induction, there is $0 \neq \overline{G}'' \subseteq \oE''$ saturated with $\xi_{\oG'', \oE''} \in \hC_{++} (X)$ and
\begin{equation}\label{nochnsternchen}
\xi^2_{\oG'',\oE''} \ge \frac{-\Delta''}{r^{''2} (r'' - 1)} > \frac{r^2}{r^{''2}} \xi^2_{\oE',\oE} \; .
\end{equation}
Clearly $\overline{G} := \ker (E \to E'' / G'') \subseteq \oE$ is saturated and from Lemma \ref{2}, the positivity of the coefficients appearing there and lemma \ref{cone} we get
  
\begin{eqnarray*}
  |\xi_{\oG,\oE}| & \ge & \frac{r' (r'' -s)}{(r' + s) r''} |\xi_{\oE',\oE}| + \frac{s}{r' + s} |\xi_{\oG'',\oE''}| \\
& \stackrel{(\ref{nochnsternchen})}{>} & \frac{r' (r'' - s)}{(r'+s) r''} |\xi_{\oE' , \oE}| + \frac{s}{r'+s} \frac{r}{r''} |\xi_{\oE',\oE}| \\
& = & \left( \frac{r' (r'' -s) + rs}{r'' (r' + s)} \right) |\xi_{\oE',\oE}| = |\xi_{\oE',\oE}| \; .
\end{eqnarray*}
This again contradicts the maximality of $|\xi_{\oE',\oE}|$ and concludes the proof.
\end{proofof}

\end{document}